\documentclass[10pt,a4paper]{amsart}
\usepackage{amssymb}
\usepackage{amsmath}
\usepackage{amsfonts}

\newcommand{\Z}{\mathbb{Z}}
\newcommand{\F}{\mathbb{F}}

\newtheorem{main-dummy}{Main-Dummy}
\newtheorem{dummy}{Dummy}

\newtheorem{main-theorem}[main-dummy]{Theorem}

\newtheorem{lemma}[dummy]{Lemma}
\newtheorem{theorem}[dummy]{Theorem}
\newtheorem{prop}[dummy]{Proposition}
\newtheorem{cor}[dummy]{Corollary}

\theoremstyle{definition}
\newtheorem{definition}[dummy]{Definition}

\newtheorem{example}[dummy]{Example}

\theoremstyle{remark}
\newtheorem{rem}[dummy]{Remark}

\begin{document}
\bibliographystyle{amsalpha}
\author{Sandro Mattarei}

\email{mattarei@science.unitn.it}

\urladdr{http://www-math.science.unitn.it/\~{ }mattarei/}

\address{Dipartimento di Matematica\\
  Universit\`a degli Studi di Trento\\
  via Sommarive 14\\
  I-38050 Povo (Trento)\\
  Italy}

\title[Linear recurrence relations for modular binomial coefficients]
{Linear recurrence relations for binomial coefficients modulo a prime}

\begin{abstract}
We investigate when the sequence of binomial coefficients $\binom{k}{i}$ modulo a prime $p$,
for a fixed positive integer $k$,
satisfies a linear recurrence relation of (positive) degree $h$ in the finite range $0\le i\le k$.
In particular, we prove that this cannot occur if $2h\le k<p-h$.
This hypothesis can be weakened to $2h\le k<p$ if we assume, in addition, that the characteristic polynomial
of the relation does not have $-1$ as a root.
We apply our results to recover a known bound for the number of points of a Fermat curve over a finite field.
\end{abstract}

\date{10 August 2005}



\thanks{The  author  is grateful  to  Ministero dell'Istruzione, dell'Universit\`a  e
  della  Ricerca, Italy,  for  financial  support of the
  project ``Graded Lie algebras  and pro-$p$-groups of finite width''.}

\maketitle

\thispagestyle{empty}

\section{Introduction}\label{sec:intro}

As is customary, let the binomial coefficients $\binom{k}{i}$ be defined by the identity
\[
(1+x)^k
=\sum_{i\in\Z}
\binom{k}{i}x^i
=\sum_{i\ge 0}
\binom{k}{i}x^i
\]
in the ring of formal power series $\Z[[x]]$, for $k,i\in\Z$.
In particular, $\binom{k}{i}$ vanishes if $i<0$ or if $i>k$.
Consider the sequence $\binom{k}{i}$ for a fixed $k$.
It is clearly never periodic on the whole range $i\in\Z$, and
restricted to the range $i\ge 0$ it is periodic exactly
in one case, namely $k=-1$, where $\binom{k}{i}=(-1)^i$.
Replicas of this isolated instance in characteristic zero appear when the sequence $\binom{k}{i}$
is viewed modulo a prime $p$:
the reduced sequence is periodic with period two in the natural range $0\le i\le k$ if $k+1$ is a power of $p$.
In fact, because of the identity $(a+b)^p=a^p+b^p$ in characteristic $p$, in the ring $\F_p[x]$ we have
\[
(1+x)^{p^s-1}=
(1+x^{p^s})/(1+x)
=(1+x^{p^s})\sum_{i\ge 0}
(-1)^ix^i
=\sum_{i=0}^{p^s-1}
(-1)^ix^i
\]
and hence $\binom{p^s-1}{i}\equiv (-1)^i\pmod{p}$ for $0\le i\le k$.
We have proved in~\cite{Mat:binomial} that under some fairly natural further assumptions this is
the only occurrence of periodicity for the sequence $\binom{k}{i}$ modulo $p$ in the range ${0\le i\le k}$,
for a fixed $k\ge 0$.
In particular, Corollary~4.2 of~\cite{Mat:binomial} asserts that
if $k+1$ is not a power of $p$ then
the sequence of binomial coefficients $\binom{k}{i}$ modulo $p$, considered in the range ${0\le i\le k}$,
cannot be periodic of any period $h$ prime to $p$ and with $2h\le k$.
A similar assertion holds for the signed binomial coefficients $(-1)^i\binom{k}{i}$.
In fact, both assertions hold under weaker and more precise assumptions,
for which we refer to~\cite{Mat:binomial}.

Since a periodicity relation is the special case of a linear recurrence relation
where the characteristic polynomial has the form $x^h-1$, it is natural to
ask when the sequence $\binom{k}{i}$ modulo $p$, for $k$ fixed,
satisfies a linear recurrence relation in the range $0\le i\le k$.
Note that, in characteristic zero and in the range $i\ge 0$, the sequence $\binom{k}{i}$
satisfies the linear recursion with characteristic polynomial $(1+x)^{-k}$ when $k<0$
(see Example~\ref{ex:k_large} for a similar instance with $k\ge 0$ in positive characteristic),
and the linear recursion with characteristic polynomial $x^{k+1}$ for $k\ge 0$
(because the sequence vanishes for $i>k$).
However, the problem becomes more interesting when we restrict our attention to the finite range $0\le i\le k$,
for $k\ge 0$.
Naturally, the familiar definitions pertaining to linear recurrence relations for infinite sequences
need to be adjusted to the case of finite sequences, as we do in Section~\ref{sec:recurrence}.
In particular, it will appear that a natural requirement to avoid degenerate cases
is to consider linear recurrence relations of order $h$
only for sequences of at least $2h+1$ terms, see Remark~\ref{rem:h_large} and Example~\ref{ex:h_large}.
In the case of the sequence $\binom{k}{i}$ modulo $p$ in the natural range $0\le i\le k$,
this assumption becomes $2h\le k$, which we have already encountered in Corollary~4.2 of~\cite{Mat:binomial}
quoted above.
An analogue of that result for linear recurrence instead of periodicity
is our Theorem~\ref{thm:recurrence}.
A simplified version of that asserts that
the sequence of binomial coefficients $\binom{k}{i}$ modulo $p$ restricted to the range ${0\le i\le k}$
cannot satisfy a linear recurrence relation of degree $h$ if $0<2h\le k<p-h$.
The assumption $k<p-h$, which is indispensable according to Example~\ref{ex:k_large},
can be weakened to $k<p$ provided we assume that the characteristic polynomial of the linear recurrence relation
does not have $-1$ as a root, as in Theorem~\ref{thm:recurrence_stronger}.

We base our proofs of Theorems~\ref{thm:recurrence} and~\ref{thm:recurrence_stronger}
on two different methods.
Both methods would actually work in both cases, as we explain in Remark~\ref{rem:other_proofs},
but each method may have its own strengths in view of possible generalizations,
notably to values of $k$ larger than $p$.
The first method comes naturally from the ordinary theory of linear recurrent sequences
and consists in evaluating certain {\em Hankel determinants}.
We do that in Proposition~\ref{prop:determinant}, which may be of independent interest.
The second method employs the generating function $(1+x)^k$ for the binomial coefficients
in a more explicit way.
It is based on Lemma~\ref{lemma:HBK}, a slight extension of an elementary fact taken from \cite{H-BK},
asserting that the number of nonzero coefficients of a polynomial in characteristic $p$
exceeds the multiplicity of one of its roots, provided that multiplicity is less than $p$.
The characteristic zero analogue of this fact appears as Lemma~1 in~\cite{Bri}, but may be well known.

The same arguments employed to prove Theorems~\ref{thm:recurrence} and~\ref{thm:recurrence_stronger}
can be used to recover a known bound for the number of points of a Fermat curve on a finite field.
We do this in Section~\ref{sec:Fermat}, to which we refer for an introduction to the problem.
Since this topic does not require an understanding of linear recurrence relations,
we have kept Section~\ref{sec:Fermat} essentially independent from the rest of the paper.
However, we do explain the connection with linear recurrence relations for binomial coefficients
in Remark~\ref{rem:connection}.

\section{Linear recurrence relations over a finite range}\label{sec:recurrence}

In this section we recall from~\cite{LN:Introduction} some basic concepts concerning (homogeneous) linear recurrences
satisfied by an infinite sequence, and adapt them to finite sequences.

Let $a(x)=\sum_{i\ge 0}a_ix^i$ be a monic polynomial of degree $h$ (possibly zero).
A sequence $\{s_i\}_{i\ge 0}$ satisfies the (homogeneous) linear recurrence relation
with {\em characteristic polynomial} $a(x)$ if
\begin{equation}\label{eq:recurrence}
s_{i+h}=-a_{h-1}s_{i+h-1}-\cdots-a_0s_{i}
\quad\text{for all $i\ge 0$}.
\end{equation}
The degree $h$ of $a(x)$ is called the {\em order} of the linear recurrence relation,
in analogy with the terminology used for linear differential equations.
This definition of a linear recurrence relation is motivated by applications where the recurrence allows one
to compute $s_i$ from the $h$ elements preceding it in the sequence.
However, the {\em reciprocal characteristic polynomial} $a^\ast(x)=x^ha(1/x)=\sum_{i\ge 0}a_{h-i}x^i$
is somehow more suited to algebraic manipulations than the characteristic polynomial.
(Note that $a^\ast(x)$ may have degree lower than $a(x)$, namely, when $a_0=0$.)
In particular, setting $a_i^\ast=a_{h-i}$, and hence
$a^\ast(x)=\sum_{i\ge 0}a_i^\ast x^i$,
we can rewrite~\eqref{eq:recurrence} as
\begin{equation*}
a_h^\ast s_{i}+a_{h-1}^\ast s_{i+1}+\cdots+a_0^\ast s_{i+h}=0
\quad\text{for all $i\ge 0$}.
\end{equation*}
These equations impose the vanishing of the coefficient of $x^{i+h}$ in the product
$a^\ast(x)s(x)$, for all $i\ge 0$,
where $s(x)=\sum_{i\ge 0}s_ix^i$ is the {\em generating function} of the sequence $\{s_i\}_{i\ge 0}$.
Therefore, a sequence $\{s_i\}_{i\ge 0}$ satisfies the linear recurrence associated with $a(x)$ if and only if
$a^\ast(x)s(x)$ is a polynomial of degree less than $h$
(cf.~\cite[Theorem~6.40]{LN:Introduction}).

Consider now a finite sequence $\{s_i\}_{u\le i\le v}$, and let
$s(x)=\sum_{u\le i\le v}s_ix^i$
be its generating function.
We may encompass the classical case of infinite sequences by allowing $v=\infty$.
Since we are not assuming that $u\ge 0$, in general $s(x)$
is a formal Laurent series rather than an ordinary power series.

\begin{definition}\label{def:recurrence}
Let $a(x)=\sum_{i\ge 0}a_ix^i$ be a monic polynomial of degree $h$.
The sequence $\{s_i\}_{u\le i\le v}$
satisfies the linear recurrence with {\em characteristic polynomial}  $a(x)$ if
\begin{equation}\label{eq:recurrence_finite}
a_0s_{i}+a_1s_{i+1}+\cdots+a_hs_{i+h}=0
\quad\text{for $u\le i\le v-h$}.
\end{equation}
\end{definition}

A finite sequence $\{s_i\}_{u\le i\le v}$ satisfies
the linear recurrence with characteristic polynomial $a(x)$
if and only if the sequence can be extended to an infinite sequence $\{s_i\}_{i\ge u}$
satisfying the linear recurrence in the usual sense.
This reveals a slight asymmetry in Definition~\ref{def:recurrence} with respect to reversing the ordering
of the finite sequence, due to our requirement that $a_h$ be nonzero (and then, without loss, being equal to one)
without a similar requirement on $a_0$.
We accept to live with this harmless asymmetry rather than departing from the standard terminology
used for infinite sequences.

Let
$a^\ast(x)=x^ha(1/x)=\sum_{i\ge 0}a_{h-i}x^i$
be the {\em reciprocal characteristic polynomial}.
Arguing as in an earlier paragraph
we find that~\eqref{eq:recurrence_finite} is satisfied
if and only if the coefficient of $x^j$ in the polynomial $a^\ast(x)s(x)$ vanishes for
$u+h\le j\le v$.
In particular, a finite sequence $\{s_i\}_{u\le i\le v}$ satisfies vacuously any linear recurrence relation
of order $h>v-u$.
It also follows easily that if the sequence satisfies a linear recurrence relation with characteristic polynomial $a(x)$
then it satisfies any linear recurrence relation whose characteristic polynomial is a multiple of $a(x)$.

The {\em Hankel determinants}
$D_r^{(h)}=\det\big((s_{r+i+j})_{i,j=0,\ldots,h-1}\big)$
play an important role in the ordinary theory of infinite linear recurring sequences,
see~\cite[Chapter~6]{LN:Introduction}.
With some care some of their properties can be translated to the present setting of finite sequences.
Here we limit ourselves to the following basic fact.

\begin{lemma}\label{lemma:Hankel}
If the sequence $\{s_i\}_{u\le i\le v}$
satisfies a linear recurrence of order $h$ then
$D_r^{(h+1)}=0$ for $u\le r\le v-2h$.
\end{lemma}

\begin{proof}
View~\eqref{eq:recurrence_finite} as a system of $v-h-u+1$ homogeneous linear equations
in the $h+1$ indeterminates $a_0,\ldots,a_h$.
The Hankel determinants $D_r^{(h+1)}$ under consideration are the determinants of the subsystems consisting
of $h+1$ consecutive equations.
The existence of a nonzero solution (with $a_h=1$ here) implies that the matrix of the system has rank less than $h+1$,
and hence all Hankel determinants vanish.
\end{proof}

Note that the values for $r$ in Lemma~\ref{lemma:Hankel} are all those for which $D_r^{(h+1)}$ is defined.
In particular, the conclusion of Lemma~\ref{lemma:Hankel} is void if $2h>v-u$.

\begin{rem}\label{rem:h_large}
Given a finite sequence $\{s_i\}_{u\le i\le v}$ and a positive integer $h$
with $h\le v-u<2h$, the set of $v-h-u+1$ equations given by~\eqref{eq:recurrence_finite} necessarily
has a nonzero solution $(a_0,\ldots,a_h)$, but need not have any with $a_h\neq 0$.
In fact, the sequence need not satisfy any linear recurrence relation of order $h$ in this case,
as is shown by the sequence with $s_{v}=1$ and $s_i=0$ for $u\le i<v$.
Nevertheless, the condition $2h\le v-u$ is a natural assumption when claiming that a finite sequence
{\em does not} satisfy a linear recurrence relation, as in Theorems~\ref{thm:recurrence}
and~\ref{thm:recurrence_stronger} below.
We examine a specific instance in Example~\ref{ex:h_large}.
\end{rem}

\section{Linear recurrence relations for binomial coefficients}\label{sec:binomial}

One way of studying linear recurrence relations for the sequence of binomial coefficients $\binom{k}{i}$
(for $k$ fixed) is evaluating the corresponding Hankel determinants.
We use the notation $k^{\underline{r}}=k(k-1)\cdots(k-r+1)$ for $k$ and $r$ integers with $r\ge 0$,
reading $k^{\underline{0}}=1$.

\begin{prop}\label{prop:determinant}
Let $k$ be an integer and let $h,r$ be nonnegative integers.
Let
$B(k,h,r)$ denote the matrix $\left(\binom{k}{r+i+j}\right)_{i,j=0,\ldots,h}$.
Then we have
\begin{equation*}
\det(B(k,h,r))=
(-1)^{\binom{h+1}{2}}
\prod_{s=0}^h
\frac{(k+s)^{\underline{r+h}}}
{(r+h+s)^{\underline{r+h}}}.
\end{equation*}
\end{prop}

\begin{proof}
We proceed by induction on $h$.
The conclusion holds for $h=0$ since $\binom{k}{r}=\frac{k^{\underline{r}}}{r^{\underline{r}}}$.
Now assume $h>0$.
To avoid confusion, we count the rows and columns of a matrix according to their index.
Thus, we call $0$th row the earliest row of $B(k,h,r)$.
We compute the determinant according to Laplace's rule, with respect to the last column (that is, the $h$th column).
For $i=h, h-1,\ldots,1$ (in this order) we subtract from the $i$th row the $(h-1)$st row multiplied by
$\frac{k+1}{r+h+i}-1$.
Using the binomial identities
$\binom{a}{b+1}+\binom{a}{b}=\binom{a+1}{b+1}$
and
$\binom{a}{b}\frac{a+1}{b+1}=\binom{a+1}{b+1}$
we find that
the $(i,j)$-entry of the resulting matrix, for $i>0$, equals
\begin{align*}
&\binom{k}{r+i+j}-\binom{k}{r+i+j-1}
\left(\frac{k+1}{r+h+i}-1\right)
\\
&\qquad=
\binom{k+1}{r+i+j}
-\binom{k}{r+i+j-1}
\frac{k+1}{r+h+i}
\\
&\qquad=
\binom{k+1}{r+i+j}
\left(1-
\frac{r+i+j}{r+h+i}
\right)
=
\binom{k+1}{r+i+j}
\frac{h-j}{r+h+i}.
\end{align*}
In particular, the $h$-th column of $B(k,h,r)$ vanishes except for its $(0,h)$-entry, which equals
$\binom{k}{r+h}$.
Consequently, $\det(B(k,h,r))$ equals $(-1)^h\binom{k}{r+h}$
times the determinant of the matrix
obtained by removing from it the $0$th root and the $h$th column.
After shifting its row-index by one, the latter matrix becomes
\begin{equation}\label{eq:matrix}
\left(\binom{k+1}{r+1+i+j}\frac{h-j}{r+h+i+1}\right)_{i,j=0,\ldots,h-1}.
\end{equation}
By collecting the factor $1/(r+h+i+1)$ from the $i$th row
and the factor $h-j$ from the $j$th column, for each row and column, we find that the determinant of the matrix
in~\eqref{eq:matrix} equals the product of
$\binom{r+2h}{h}^{-1}$ and $\det(B(k+1,h-1,r+1))$.
Since
$\binom{k}{r+h}\binom{r+2h}{h}^{-1}=
\frac{k^{\underline{r+h}}}
{(r+2h)^{\underline{r+h}}}$
we conclude that
\begin{align*}
\det(B(k,h,r))=
(-1)^h\frac{k^{\underline{r+h}}}
{(r+2h)^{\underline{r+h}}}
\cdot\det(B(k+1,h-1,r+1)).
\end{align*}
By induction hypothesis we have
\begin{align*}
\det(B(k,h,r))
&=
(-1)^{h+\binom{h}{2}}
\frac{k^{\underline{r+h}}}
{(r+2h)^{\underline{r+h}}}
\prod_{s=0}^{h-1}
\frac{(k+1+s)^{\underline{r+h}}}
{(r+h+s)^{\underline{r+h}}}
\\&=
(-1)^{h+\binom{h}{2}}
\prod_{s=0}^{h}
\frac{(k+s)^{\underline{r+h}}}
{(r+h+s)^{\underline{r+h}}},
\end{align*}
which concludes the proof.
\end{proof}

\begin{cor}\label{cor:determinant}
Let $p$ be a prime and let $k,h,r$ be integers with $k,h\ge 0$ and $0\le r+h\le k<p-h$.
Then $\det(B(k,h,r))$
is not a multiple of $p$.
\end{cor}

\begin{proof}
The conclusion follows from Proposition~\ref{prop:determinant} if $r\ge 0$.
The case where $r<0$ is reduced to the other case by means of the identity
\[
\det(B(k,h,r))=\det(B(k,h,k-r-2h)),
\]
which follows from the identity
$\binom{k}{i}=\binom{k}{k-i}$
for the binomial coefficients.
\end{proof}

The necessity of the conditions
$0\le r+h\le k$
in Corollary~\ref{cor:determinant} can also be seen by noting that the matrix
$B(k,h,-h-1)$ (respectively~$B(k,h,k-h+1)$)
has zeroes on and above (respectively~below) its secondary diagonal.

Our first main result states that the sequence of binomial coefficients
$\binom{k}{i}$, considered in the range $i=0,\ldots,k$ and reduced modulo a prime $p$,
does not satisfy any recurrence relation of degree $h$ if $2h\le k<p-h$.
However, in view of an application in the next section it is convenient to allow a more general range for $i$.

\begin{theorem}\label{thm:recurrence}
Let $p$ be a prime, $k,h$ nonnegative integers, and $u,v$ integers with
\[
-h\le u\le v\le k+h,\quad
2h\le v-u,\quad\text{and}\quad
k<p-h.
\]
Then the sequence of binomial coefficients $\binom{k}{i}$ modulo $p$, considered in the range $u\le i\le v$,
does not satisfy any linear recurrence relation of order $h$.
\end{theorem}

\begin{proof}
Suppose for a contradiction that the sequence under consideration satisfies a linear recurrence relation of order $h$.
According to Lemma~\ref{lemma:Hankel}, the Hankel determinants $D_r^{(h+1)}$ of the sequence vanish for $u\le r\le v-2h$.
Since $D_r^{(h+1)}$ is the reduction modulo $p$
of $\det(B(k,h,r))$, it follows in particular that
$p$ divides $\det(B(k,h,u))$.
This contradicts Corollary~\ref{cor:determinant}.
\end{proof}

The extreme case $h=0$ of Theorem~\ref{thm:recurrence} amounts to the simple fact
that $p$ does not divide $\binom{k}{i}$ for $0\le i\le k<p$.
Corollary~\ref{cor:determinant} and Theorem~\ref{thm:recurrence}
can be interpreted in characteristic zero by reading $p=\infty$
(or, equivalently, by disregarding the hypotheses which involve $p$).

We present a couple of examples to justify the hypotheses $2h+u+v\le k$ and $k<p-h$
in Theorem~\ref{thm:recurrence}, which are of a quite different nature.

\begin{example}\label{ex:h_large}
Work in characteristic zero first, and
let $k,h$ be integers with $0<2h-1=k$.
Then the sequence of binomial coefficients $\binom{k}{i}$, considered in the range ${0\le i\le k}$,
satisfies a unique linear recurrence relation of order $h$.
In fact, we may set $a_h=1$ and view~\eqref{eq:recurrence_finite} as a system of $k+1=2h$ linear equations
in the $h$ indeterminates $a_0,\ldots,a_{h-1}$.
Since its matrix $B(k,h-1,0)$ is nonsingular according to
Proposition~\ref{prop:determinant}, the system has a unique solution,
which can be computed by means of the
Berlekamp-Massey algorithm described in~\cite[Chapter~6, \S6]{LN:Introduction}.
Consequently, the reduction of the sequence modulo any prime $p$ also satisfies a linear recurrence, and this is unique
if $p\ge k+h$ according to Corollary~\ref{cor:determinant}.
\end{example}

\begin{example}\label{ex:k_large}
If $q$ is a power of the prime $p$ and $0<h\le q$ then the sequence of binomial coefficients $\binom{q-h}{i}$
modulo $p$, considered in the range $-h+1\le i\le q-1$
(which includes the more natural range $0\le i\le q-h$), satisfies the linear recurrence relation
with (reciprocal) characteristic polynomial $(1+x)^h$.
This is so because the product of its generating function $(1+x)^{q-h}$ with $(1+x)^h$
equals $1+x^q$, whose coefficients of degrees $1,\ldots,q-1$ vanish.
\end{example}

As Example~\ref{ex:k_large} suggests, it is possible to weaken the hypothesis
that $k<p-h$ in Theorem~\ref{thm:recurrence} to $k<p$ provided one assumes that the characteristic
polynomial of the linear recurrence relation does not have $1$ as root.
This can be established by a variation of the method of proof of Theorem~\ref{thm:recurrence},
which we sketch in Remark~\ref{rem:other_proofs} below.
However, it is simpler to base a proof on a different method.
We need the following refinement of Lemma~6 of~\cite{H-BK}, which had the stronger hypothesis
$\deg(f)<p$.
We call {\em weight} of a polynomial the number of its nonzero coefficients.

\begin{lemma}\label{lemma:HBK}
Let $f(x)$ be a polynomial over a field of characteristic $p$
having a nonzero root $\xi$ with multiplicity exactly $k$, with $0<k<p$.
Then $f(x)$ has weight at least $k+1$.
\end{lemma}

\begin{proof}
If $\xi$ is a root of $f(x)$ with multiplicity exactly $k$, then $f(\xi^{-1}x)$
has $1$ as a root with the same multiplicity, and has the same weight as $f(x)$.
Hence we may assume that $\xi=1$.
We proceed by induction on $k$.
The case $k=1$ being obvious, assume that $k>1$.
By dividing $f(x)=\sum_i f_ix^{i}$ by a suitable power of $x$, which leaves its weight unchanged,
we may assume that $f_0\neq 0$.
Since $p$ does not divide $k$, the derivative
\[
f'(x)=\sum_i if_ix^{i-1}
\]
has $1$ as a root with multiplicity exactly $k-1$,
and has weight one less than the weight of $f(x)$.
By induction, $f'(x)$ has weight at least $k$, and hence
$f(x)$ has weight at least $k+1$.
\end{proof}

\begin{theorem}\label{thm:recurrence_stronger}
Let $p$ be a prime, $k,h$ nonnegative integers, and $u,v$ integers with
\[
-h\le u\le v\le k+h,\quad
2h\le v-u,\quad\text{and}\quad
k<p.
\]
Then the sequence of binomial coefficients $\binom{k}{i}$ modulo $p$, considered in the range $u\le i\le v$,
does not satisfy any linear recurrence relation of order $h$ with characteristic polynomial prime to $x+1$.
\end{theorem}

\begin{proof}
Suppose that the sequence under consideration satisfies a linear recurrence relation
with characteristic polynomial $a(x)$, of degree $h$ and not having $-1$ as a root.
Therefore, the coefficient of $x^j$ in the polynomial $a^\ast(x)(1+x)^k$ vanishes for
$u+h\le j\le v$.
Consequently, $a^\ast(x)(1+x)^k$ has weight at most
$k+2h-v+u$.
However, according to Lemma~\ref{lemma:HBK} the weight of
$a^\ast(x)(1+x)^k$ is at least $k+1$.
It follows that
$2h\ge v-u+1$, which contradicts one of our hypotheses.
\end{proof}

\begin{rem}\label{rem:other_proofs}
We mentioned above that the method of proof of Theorem~\ref{thm:recurrence}
can be adapted to give a proof of Theorem~\ref{thm:recurrence_stronger}.
We briefly sketch the corresponding argument.
Assume that $p-h\le k<p$, otherwise Theorem~\ref{thm:recurrence} applies.
As in the proof of Lemma~\ref{lemma:Hankel}
one may view~\eqref{eq:recurrence_finite} as a system of $v-h-u+1$ homogeneous linear equations
in the $h+1$ indeterminates $a_0,\ldots,a_h$.
The matrix of the subsystem formed by the first $h+1$ equations equals the reduction modulo $p$ of $B(k,h,u)$,
whose determinant vanishes according to Corollary~\ref{cor:determinant}.
However, one can use Corollary~\ref{cor:determinant} to show that
the matrix has rank at least $p-k$.
(It will turn out that the matrix has rank exactly $p-k$.)
Therefore, the space of solutions of the system~\eqref{eq:recurrence_finite}
has dimension at least $h+1-p+k$.
According to Example~\ref{ex:k_large} and an earlier observation,
the sequence under consideration satisfies any linear recurrence relation which has a multiple of
$(1+x)^{p-k}$ as characteristic polynomial.
Those characteristic polynomials (assumed monic here) which have degree $h$
form an affine subspace of $\F_p[x]$ of dimension $h-p+k$, and hence
span a linear subspace of dimension $h-p+k+1$.
It follows that these account for all solutions of the system~\eqref{eq:recurrence_finite},
which is the desired conclusion.
\end{rem}

It is also possible to use Lemma~\ref{lemma:HBK} to prove Theorem~\ref{thm:recurrence}.
In fact, the weaker version of Lemma~\ref{lemma:HBK} given in~\cite{H-BK}, where $\deg(f)<p$, would suffice for that.

\section{An application to Fermat curves over a finite field}\label{sec:Fermat}

Let $\F_q$ be the finite field of $q$ elements  and let $p$ be its characteristic.
Consider the {\em Fermat curve}
$ax^n+by^n=z^n$, expressed in homogeneous coordinates, where $n>1$ is an integer prime to $p$,
and $a,b\in\F_q^\ast=\F_q\setminus\{0\}$.
A classical estimate on the number $N_n(a,b,q)$ of its projective $\F_q$-rational points is
\[
|N_n(a,b,q)-q-1|\le (n-1)(n-2)\sqrt{q}.
\]
This is originally due to Hasse and Davenport~\cite{DavHas} but is a special case of
{\em Weil's bound} for curves over finite fields.
Weil's bound for Fermat curves is easy to prove by means of Gauss and Jacobi sums,
as well as its generalisation to {\em diagonal equations} in several variables,
see~\cite{IR}, \cite{LN} or~\cite{Small}.
An alternative proof is based on the character theory of a finite Frobenius group,
see~\cite[Section~26]{Feit} for the basic argument and~\cite{Mat:Fermat-character}
for a refinement.

Weil's upper bound for $N_n(a,b,q)$ is not optimal when $n$ (and with it the genus of the curve)
is relatively large with respect to $q$.
Better upper bounds in this situation were found by Garc\'{\i}a and Voloch, using tools from algebraic geometry.
According to~\cite[Corollary~1]{GarVol}, rewritten here after elementary calculations,
if $s$ is an integer such that $1\le s\le n-3$ and $sn\le p$ then
\begin{equation}\label{eq:GV-bounds}
N_n(a,b,q)\le
\left(\frac{s^2-s-2}{4}+\frac{4}{s+3}
\right)n^2
+2\frac{n(q-1-d)}{s+3}+d,
\end{equation}
where $d$ is the number of $\F_q$-rational points of the curve with $xyz=0$.
Garc\'{\i}a and Voloch pointed out that
their bounds~\eqref{eq:GV-bounds} hold in more general circumstances where
the assumption $sn\le p$ may not be satisfied,
and described those circumstances in detail for the cases $s=1,2$.
In particular, the case $s=1$ of~\eqref{eq:GV-bounds}, which reads
\begin{equation}\label{eq:GV-bound}
N_n(a,b,q)\le(n(n+q-1)-d(n-2))/2,
\end{equation}
is valid without the assumption $n\le p$, for $p$ odd, except when $n$ has the form
$n=(q-1)/(r-1)$ for some subfield $\F_r$ of $\F_q$ (which are true exceptions).
Bound~\eqref{eq:GV-bound} is better than Weil's upper bound, roughly, when $n$ is larger than $\sqrt{q}/2$.
In Corollary~\ref{cor:Fermat} we establish bound~\eqref{eq:GV-bound} under the assumptions that
$n$ divides $q-1$ (which is harmless in view of the next paragraph) and that $n>(q-1)/(p-1)$.

The set $G$ of $n$th powers in $\F_q^\ast$ coincides with the set of $m$th powers, where $m=(n,q-1)$,
and is the subgroup of $\F_q^\ast$ of order $(q-1)/m$.
Setting
$(\alpha,\beta)\mapsto a\alpha^n$
gives an $m^2$-to-one map of the set of pairs $(\alpha,\beta)\in\F_q^\ast\times\F_q^\ast$
with $a\alpha^n+b\beta^n=1$ onto the set $aG\cap(1-bG)$.
Consequently, we have
$N_n(a,b,q)=m^2|aG\cap(1-bG)|+d$,
where $d$ is the number of projective $\F_q$-rational points of the curve $ax^n+by^n=z^n$ with $xyz=0$.
However, it easy to see that $d$ coincides with the number of projective $\F_q$-rational points
of the curve $ax^m+by^m=z^m$ with $xyz=0$,
and hence $N_n(a,b,q)=N_m(a,b,q)$.
Since Garc\'{\i}a and Voloch's bounds~\eqref{eq:GV-bounds} (as well as Weil's bound) do not increase
by replacing $n$ with $m$ (and leaving $d$ unchanged), it is no loss
to assume that $n$ divides $q-1$ in the sequel.
In the next result we use Lemma~\ref{lemma:HBK} to produce an upper bound for $|aG\cap(1-bG)|$.
The first part of the argument is analogous to the proof of Theorem~2 in~\cite{Mat:binomial}.

\begin{theorem}\label{thm:subgroup}
Let $G$ be a subgroup of $\F_q^\ast$ with $|G|<p-1$, and let $a,b\in\F_q^\ast$.
Set $e=0,1,2,3$ according as none, one, two or all three of $a$, $b$ and $-a/b$ (counting repetitions) belong to $G$.
Then $|a G\cap(1-b G)|\le (|G|+1-e)/2$.
\end{theorem}

\begin{proof}
The elements of the cosets $a G$ and $b G$ of $G$ in $\F_q^\ast$ are the roots of the polynomials
$x^k-a^k$ and $x^k-b^k$, where $k=|G|$.
Consequently, the elements of $a G\cap(1-b G)$ are the roots of the greatest common divisor
$(x^k-a^k,(1-x)^k-b^k)$, which we write in the form
$(x^k-a^k)/f(x)$, where $f(x)=\prod_{\xi\in a G\setminus(1-b G)}(x-\xi)$.
Hence $|a G\cap(1-b G)|=k-h$, where $h=\deg(f)$.
There exists a polynomial $g(x)\in\F_q[x]$, necessarily of degree $h$ and with leading coefficient $(-1)^k$,
such that
\begin{equation*}
(1-x)^k-b^k=
g(x)
\frac{x^{k}-a^k}{f(x)},
\end{equation*}
and hence
\begin{equation}\label{eq:expanded}
f(x)(1-x)^k
=b^kf(x)-a^kg(x)+x^{k}g(x).
\end{equation}
The polynomial $f(x)(1-x)^k$ has $1$ as a root with multiplicity exactly $k$ or $k+1$ according as
$a\not\in G$ or $a\in G$.
According to Lemma~\ref{lemma:HBK}, its weight is at least $k+1$ in the former case, and at least $k+2$ in the latter.
However, the polynomial at the right-hand side of~\eqref{eq:expanded} has weight at most $2h+2$.
Consequently, in any case we have $2h+2\ge k+1$, that is, $k-h\le(k+1)/2$.
This is the desired conclusion in case $e=0$.

The remaining cases are established by taking into account whether $a\in G$,
and noting that the right-hand side of~\eqref{eq:expanded} has actually weight at most $2h+1$
if either $b$ or $-a/b$ belongs to $G$, and at most $2h$ if both do.
In fact, if $b\in G$ then $g(x)$ has $0$ as a root, and hence has no constant term,
while if $-a/b\in G$ then $b^k f(x)$ and $a^k g(x)$ have the same leading coefficient.
\end{proof}

\begin{rem}\label{rem:connection}
We sketch a minor variation of the proof of Theorem~\ref{thm:subgroup} which emphasizes the connection
with the linear recurrence relations for binomial coefficients discussed in the previous section.
For simplicity we restrict ourselves to the case $e=0$.
Expanding the product on the left-hand side of~\eqref{eq:expanded} and writing $f(-x)=\sum_{j=0}^h f_jx^j$
we obtain that
\begin{equation}\label{eq:dependence}
f_h\binom{k}{s}+f_{h-1}\binom{k}{s+1}+\cdots+f_0\binom{k}{s+h}=0
\end{equation}
for each integer $s$ such that $x^{s+h}$ has coefficient zero in the polynomial at the right-hand side
of~\eqref{eq:expanded}.
This certainly holds for $1\le s\le k-h-1$, and hence the sequence of binomial coefficients $\binom{k}{i}$
modulo $p$, restricted to the range $1\le i\le k-1$,
satisfies a linear recurrence relation of order $h$, with characteristic polynomial prime to $x+1$.
According to Theorem~\ref{thm:recurrence_stronger} we have $2h>k-2$, and the conclusion follows.
\end{rem}

\begin{cor}\label{cor:Fermat}
Let $q$ be a power of the prime $p$, $n$ a divisor of $q-1$ with $n>(q-1)/(p-1)$, and $a,b\in\F_p^\ast$.
Then the Fermat curve $ax^n+by^n=z^n$ has at most
$(n(n+p-1)-d(n-2))/2$ projective $\F_p$-rational points,
where $d$ is the number of points with $xyz=0$.
\end{cor}

\begin{proof}
The number of projective $\F_q$-rational points of the curve with $xyz\neq 0$
equals $n^2|aG\cap(1-bG)|$, where $G$ is the subgroup of $\F_q^\ast$ of order $(q-1)/n$.
According to Theorem~\ref{thm:subgroup}, this number is at most
$(n(q-1+n-en)/2$.
Adding to this the number of points with $xyz=0$, which is $d=en$, we reach the conclusion.
\end{proof}

\bibliography{References}

\end{document}